\documentclass[preprints,article,accept,oneauthor,pdftex]{my-mdpi} 

\firstpage{1} 
\pubvolume{10}
\issuenum{3}
\articlenumber{171}
\pubyear{2021}
\copyrightyear{2021}
\externaleditor{Academic Editor: Giampiero Palatucci} 
\datereceived{12 June 2021} 
\daterevised{22 July 2021} 
\dateaccepted{26 July 2021} 
\datepublished{29 July 2021} 
\hreflink{https://doi.org/10.3390/\newline axioms10030171} 
\doinum{10.3390/axioms10030171}

\pdfoutput=1

\usepackage{MnSymbol} 
\usepackage[mathscr]{euscript}

\Title{On a Non-Newtonian Calculus of Variations}

\TitleCitation{On a Non-Newtonian Calculus of Variations}

\Author{Delfim F. M. Torres \orcidA{}}

\AuthorNames{Delfim F. M. Torres}

\AuthorCitation{Torres, D.F.M.}

\address[1]{Center for Research and Development in Mathematics and Applications (CIDMA),
Department of Mathematics, \mbox{University of Aveiro}, 3810-193 Aveiro, Portugal; 
delfim@ua.pt}

\abstract{The calculus of variations is a field of mathematical analysis  
born in 1687 with Newton's problem of minimal resistance, which is concerned 
with the maxima or minima of integral functionals. Finding the solution 
of such problems leads to solving the associated Euler--Lagrange equations.
The subject has found many applications over the centuries, e.g., in
physics, economics, engineering and biology. Up to this moment,
however, the theory of the calculus of variations has been confined
to Newton's approach to calculus. As in many applications
negative values of admissible functions are not physically plausible, 
we propose here to develop an alternative calculus of variations
based on the non-Newtonian approach first introduced by Grossman and Katz 
in the period between 1967 and 1970, which provides a calculus defined, 
from the very beginning, for positive real numbers only, and it is based
on a (non-Newtonian) derivative that permits one to compare relative changes 
between a dependent positive variable 
and an independent variable that is also positive.
In this way, the non-Newtonian calculus of variations
we introduce here provides a natural framework for problems involving functions
with positive images. Our main result is a first-order optimality
condition of Euler--Lagrange type. The new calculus of variations 
complements the standard one in a nontrivial/multiplicative way, guaranteeing 
that the solution remains in the physically admissible positive range.
An illustrative example is given.}



\begin{document}


\section{Introduction}

A popular method of creating a new mathematical system 
is to vary the axioms of a known one. Non-Newtonian calculi 
provide alternative approaches to the usual calculus 
of Newton (1643--1727) and Leibniz (1646--1716), which were
first introduced by Grossman and Katz (1933--2010) 
in the period between 1967 and 1970~\cite{MR0430173}.
The two most popular non-Newtonian calculi
are the multiplicative and bigeometric calculi, which
in fact are modifications of each other: in these calculi, 
the  addition and subtraction are changed 
to multiplication and division~\cite{MR4164307}.
Since such multiplicative calculi are variations
on the usual calculus, the~traditional one is sometimes
called the additive calculus~\cite{MR2864779}.

Recently, it has been shown that non-Newtonian/multiplicative calculi 
are more suitable than the ordinary Newtonian/additive calculus
for some problems, e.g.,~in actuarial science, finance, economics, 
biology, demography, pattern recognition in images,
signal processing, thermostatistics and quantum information 
theory~\cite{MR2864779,MORA20121245,MR3452941,e22101180,MR4188122}. 
This is explained by the fact that while the basis 
for the standard/additive calculus is the representation 
of a function as locally linear, the~basis of 
a multiplicative calculus is the representation 
of a function as locally exponential~\cite{MR0430173,MR2864779,MR4188122}.
In fact, the~usefulness of product integration goes back
to Volterra (1860--1940), who introduced in 1887
the notion of a product integral and used it to study solutions 
of differential equations~\cite{MR3271405,MR3992455}.
For readers not familiar with product integrals,  
we refer to the book~\cite{MR2917851},
which contains short biographical sketches of Volterra, 
Schlesinger and other mathematicians involved in the development 
of product integrals, and~an extensive list of references, 
offering a gentle opportunity to become acquainted 
with the subject of non-Newtonian integration.
For our purposes, it is enough to understand that 
a non-Newtonian calculus is a methodology that allows one 
to have a different look at problems that can be investigated 
via calculus: it provides differentiation and integration tools, 
based on multiplication instead of addition, and~in some 
cases---mainly problems of price elasticity, multiplicative growth, 
etc.---the use of such multiplicative calculi 
is preferable to the traditional Newtonian 
calculus~\cite{cordova2006multiplicative,MR2356052,MR3722924,MR4247781}.
Moreover, a~non-Newtonian calculus is a self-contained system, 
independent of any other system of calculus~\cite{MR3528851}.

The main aim of our work was to obtain, for~the first time in the literature, 
a non-Newtonian calculus of variations that involves the minimization
of a functional defined by a non-Newtonian integral with an 
integrand/Lagrangian depending on the non-Newtonian derivative. 
The calculus of variations is a field of mathematical analysis 
that uses, as~the name indicates, variations, 
which are small changes in functions, to~find maxima and minima 
of the considered functionals: mappings from a set of functions 
to the real numbers. In~the non-Newtonian framework, instead of 
the classical variations of the form $y(\cdot) + \epsilon h(\cdot)$, 
proposed by Lagrange (1736--1813) and still used nowadays in all 
recent formulations of the calculus of 
variations~\cite{MR641031,MR2004181,MR2437854}, for~example,
in the fractional calculus of variations~\cite{MR3822307,MR4159537}, 
quantum variational calculus~\cite{MR3184533,MR3771533} 
and the calculus of variations on time scales~\cite{MR2671876,MR3718404}, 
we propose here to use ``multiplicative variations''. More precisely, 
in~contrast with the calculi of variations found in the literature,
we show here, for~the first time, how to consider variations 
of the form $y(\cdot) \cdot \epsilon^{\ln{h(\cdot)}}$.
The functionals of the calculus of variations are expressed as definite integrals,
here in a non-Newtonian sense, involving functions and their derivatives,  
in~a non-Newtonian sense here. The~functions that maximize or minimize 
the functionals of the calculus of variations are found using the Euler--Lagrange 
equation, which we prove here in the non-Newtonian setting.
Given the importance of the calculus of variations in applications,
for example, in~physics~\cite{MR2323273,MR3169197}, 
economics~\cite{MR2405376,MR3313737} 
and biology~\cite{MR3157168,MR4110649},
and the importance that non-Newtonian calculus already has
in these areas, we trust that the calculus of variations  
initiated here will call attention to the research community. 
We give credit to the citation found in the 1972 book 
of Grossman and Katz~\cite{MR0430173}: 
``for each successive class of phenomena, 
a new calculus or a new~geometry''.


\section{Materials and~Methods}
\label{sec:2}

From 1967 till 1970, Grossman and Katz gave definitions 
of new kinds of derivatives and integrals, converting the 
roles of subtraction and addition into division and multiplication, 
respectively, and~established a new family of calculi, called 
non-Newtonian calculi~\cite{MR0430173,MR3609374,MR3750265},
which are akin to the classical calculus developed by
Newton and Leibniz three centuries ago.
Non-Newtonian calculi use different types of arithmetic
and their generators. Let $\alpha$ be a bijection between 
subsets $X$ and $Y$ of the set of real numbers $\mathbb{R}$, 
and endow $Y$, with~the induced operations sum and multiplication 
and the ordering given by the inverse map $\alpha^{-1}$. 
Then the $\alpha$-arithmetic is a field with the order topology
~\cite{MR3073475}. In~concrete, given a bijection 
$\alpha : X \rightarrow Y \subseteq \mathbb{R}$,
called a generator~\cite{MR3528851}, we say that $\alpha$
defines an arithmetic if the following four operations are defined:
\begin{equation}
\label{eq:1}
\begin{split}
x \oplus y &= \alpha\left(\alpha^{-1}(x)+\alpha^{-1}(y)\right),\\
x \ominus y &= \alpha\left(\alpha^{-1}(x)-\alpha^{-1}(y)\right),\\
x \odot y &= \alpha\left(\alpha^{-1}(x)\cdot\alpha^{-1}(y)\right),\\
x \oslash y &= \alpha\left(\alpha^{-1}(x)\slash\alpha^{-1}(y)\right).
\end{split}
\end{equation}

If $\alpha$ is chosen to be the identity function and $X = \mathbb{R}$,
then \eqref{eq:1} reduces to the four operators studied in school; i.e.,
one gets the standard arithmetic, from~which the traditional (Newton--Leibniz)
calculus is developed. For~other choices of $\alpha$ and $X$, we can get
an infinitude of other arithmetics from which Grossman and Katz produced
a series of non-Newton calculi, compiled in the seminal book of 1972~\cite{MR0430173}. 
Among~all such non-Newton calculi, recently great
interest has been focused on the Grossman--Katz calculus obtained
when~we fix $\alpha(x)=e^x$, $\alpha^{-1}(x)=\ln(x)$ and $X=\mathbb{R}^+$ 
for the set of real numbers strictly greater than 
zero~\cite{MR2356052,MR3528851,MR4188122,MR4247781}.
We shall concentrate here on one option,
originally called by Grossman and Katz the 
geometric/exponential/bigeometric 
calculus~\cite{MR0430173,MR557734,MR695495,MR3389882,MR3722924,MR4164307,MR4222466}, 
but from which other different terminology 
and small variations of the original calculus have
grown up in the literature, in~particular, the~multiplicative 
calculus~\mbox{\cite{MR2356052,MR2864779,MORA20121245,MR3271405,MR3389882,MR3528851,%
MR3452941,MR3661747,MR3794497,MR4016182}}, and~more recently, the~proportional 
calculus~\cite{cordova2006multiplicative,MR2724186,MR4188122,MR4247781},
which is essentially the bigeometric calculus of~\cite{MR695495}.
Here we follow closely this last approach---in particular, the~exposition
of the non-Newton calculus as found in~\cite{MR695495,MR4188122,MR4247781},
because it is appealing to scientists who seek ways to express laws
in a scale-free~form.

Throughout the text, we fix $\alpha(x)=e^x$, $\alpha^{-1}(x)=\ln(x)$, 
and $X=\mathbb{R}^+$. Then we get from \eqref{eq:1} the following operations:
\begin{equation}
\label{eq:2}
\begin{split}
x \oplus y &= x \cdot y,\\
x \ominus y &= \frac{x}{y},\\
x \odot y &= x^{\ln(y)},\\
x \oslash y &= x^{1\slash\ln(y)},
\quad y \neq 1.
\end{split}
\end{equation}

Let $a, b, c \in \mathbb{R}^{+}$. In~the non-Newtonian arithmetic
given by \eqref{eq:2}, the~following properties of the
$\odot$ operation hold (cf. Proposition~2.1 of~\cite{MR4247781}):
\begin{enumerate}
\item[(i)] $a \odot b = b \odot a$ (commutativity);
\item[(ii)] $a \odot (b \odot c) = (a \odot b)\odot c$ (associativity);
\item[(iii)] $a \odot e = a$ (Euler's/Napier's transcendent 
number $e$ is the neutral element for $\odot$);
\item[(iv)] if $a \neq 1$ and we define $a^{\{-1\}} = e \oslash a$, then
$a \odot a^{\{-1\}} = e$ (inverse element).
\end{enumerate}

We see that in this non-Newtonian algebra, $a = 1$ is the traditional ``zero''
(in the current arithmetic, $0$ represents $-\infty$). In~fact,
see Proposition~2.2 of~\cite{MR4247781}, one~has
\begin{enumerate}
\item[(v)] $b \odot a^{\{-1\}} = b \oslash a$;
\item[(vi)] $\left(a^{\{-1\}}\right)^{\{-1\}} = a$;
\item[(vii)] $\ln\left(a \odot b\right) = \ln(a) \oplus \ln(b)$;
\item[(viii)] $(a \odot b)^{\{-1\}} = a^{\{-1\}} \odot b^{\{-1\}}$.
\end{enumerate}

Based on the mentioned properties, one easily proves that
$\left(\mathbb{R}^{+},\oplus,\odot\right)$ is a field
(see \mbox{Theorem~2.3} of~\cite{MR4247781}). In~this field,
the following calculus has been developed~\cite{MR4188122,MR4247781}.

\begin{Definition}[absolute value]
\label{def:1}
The absolute value of $x\in\mathbb{R}^+$, denoted by $[[x]]$,
is given by
\begin{equation*}
[[x]] =
\begin{cases}
x & \text{if } x \geq 1,\\
1 \ominus x & \text{if } x \in (0,1).
\end{cases}
\end{equation*}
\end{Definition}

Let $x, y, z$ be positive real numbers and define
$d : \mathbb{R}^{+} \times \mathbb{R}^{+} \rightarrow \mathbb{R}^{+}$ as
$$
d(x,y) = [[x \ominus y]].
$$

The following properties are simple to~prove:
\begin{itemize}
\item $d(x,y) \geq 1$;
\item $d(x,y) = 1$ if, and~only if, $x = y$;
\item $d(x,y) = d(y,x)$;
\item $d(x,z) \leq d(x,y) + d(y,z)$.
\end{itemize}

We can now introduce the notion of~limit.

\begin{Definition}[limit]
\label{def:2}
We write 
$\lim_{x\rightarrow x_0} f(x) = L$, 
$L \in \mathbb{R}^{+}$,
as: if for all $\epsilon > 1$ there exists $\delta > 1$ such that
if $1 < d(x,x_0) < \delta$, then $d\left(f(x),L\right) < \epsilon$.
\end{Definition}

\textls[-15]{According with Definition~\ref{def:2}, it is possible to specify
the meaning of equality \linebreak \mbox{$\lim_{x\rightarrow x_0} f(x) = f(x_0)$},
and therefore, the~notion of continuity in the non-Newtonian~calculus.}

\begin{Definition}[continuity]
\label{def:3}
We say that $f$ is continuous at $x_0$ or 
$$
\lim_{x\rightarrow x_0} f(x) = f(x_0),
$$
if $\forall$ $\epsilon > 1$, $\exists$ $\delta > 1$ such that
$d(x,x_0) < \delta \implies d\left(f(x),f(x_0)\right)<\epsilon$.
\end{Definition}

\textls[-20]{{We proceed by reviewing the essentials 
on non-Newtonian differentiation \mbox{and~integration.}}}


\subsection{Derivatives}

The derivative of a function is introduced in the following~terms.

\begin{Definition}[derivative~\cite{MR695495,MR2724186,MR4188122,MR4247781}]
\label{def:4}
A positive function $f$ is differentiable at $x_0$ if
$$
\lim_{x\rightarrow x_0}\left[\left(f(x)\ominus f(x_0)\right)
\oslash\left(x\ominus x_0\right)\right]
=
\lim_{h\rightarrow 1}\left[\left(f(x_0\oplus h)\ominus f(x_0)\right)\oslash h\right]
$$
exists. In~this case, the~limit is denoted by $\widetilde{f}(x_0)$ and receives
the name of derivative of $f$ at $x_0$. Moreover, we say that $f$ is differentiable
if $f$ is differentiable at $x_0$ for all $x_0$ in the domain of $f$.
\end{Definition}

It is not difficult to prove that if $f$ is differentiable at $x_0$, 
then $f$ is continuous at $x_0$.~Define 
\begin{equation*}
\begin{split}
x^{\{0\}} &= e,\\
x^{\{n\}} &= x \odot \cdots \odot x \text{ $n$ times},
\quad n \in \mathbb{N}.
\end{split}
\end{equation*}

We have
\begin{equation}
\label{eq:ex1}
\widetilde{x^{\{n\}}} = \left(x^{\{n-1\}}\right)^n
= e^n \odot x^{\{n-1\}}, \quad n \in \mathbb{N}.
\end{equation}

In particular, if~$n = 1$ in \eqref{eq:ex1}, we~get:
\begin{itemize}
\item If $f(x) = x$, then $\widetilde{f}(x) = e$.
\end{itemize}

More examples of derivatives of a function 
in the sense of Definition~\ref{def:4} follow:
\begin{itemize}
\item If $f(x) = e^x$, then $\widetilde{f}(x) = e^x$;
\item If $f(x) = \ln(x)$, then $\widetilde{f}(x) = e \oslash x$;
\item If $\cos_e(x) := e^{\cos(\ln(x))}$ and $\sin_e(x) := e^{\sin(\ln(x))}$, 
then 
$$
\widetilde{\cos_e(x)} = 1 \ominus \sin_e(x)
\quad \text{and} \quad 
\widetilde{\sin_e(x)} = \cos_e(x).
$$
\end{itemize}

The basic rules of differentiation
(keep recalling that 1 is the ``zero''
of the non-Newtonian calculus) follow:
\begin{enumerate}
\item[(a)] If $f(x) = c$, and~$c$ is a positive constant, 
then $\widetilde{f}(x) = 1$ (derivative of a constant);
\item[(b)] $\widetilde{f\oplus g} 
= \widetilde{f} \oplus \widetilde{g}$ (derivative of a sum);
\item[(c)] $\widetilde{f\ominus g} 
= \widetilde{f} \ominus \widetilde{g}$ (derivative of a difference);
\item[(d)] $\widetilde{f\odot g} = \left(\widetilde{f} \odot g\right) 
\oplus \left(f \odot \widetilde{g}\right)$ (derivative of a product);
\item[(e)] $\widetilde{f \circ g} = \left(\widetilde{f}
\circ g\right)\odot \widetilde{g}$ (chain rule).
\end{enumerate}

If $\widetilde{f}(x) = 1$ for all $x\in (a,b)$, then $f(x) = c$ for all $x\in (a,b)$,
where $c$ is a constant. Moreover, if~$\widetilde{f}(x) = \widetilde{g}(x)$ 
for all $x\in (a,b)$, then there exists a constant $c$ such that $f(x) = c g(x)$
for all $x\in (a,b)$; that is, $f(x) = g(x) \oplus c$.

Higher-order derivatives are defined as usual:
\begin{equation*}
\begin{split}
\widetilde{f}^{(0)}(x) &= f(x),\\
\widetilde{f}^{(n)}(x) 
&= \frac{\tilde{d}}{\tilde{d} x} \left[\widetilde{f}^{(n-1)}(x)\right],
\quad n \in \mathbb{N}.
\end{split}
\end{equation*}

In the sequel, we use the following notation:
$$
{\overset{n}{\widetilde{\sum_{i=0}}}} a_i 
= a_0 \oplus \cdots \oplus a_n.
$$

\begin{Theorem}[Taylor's theorem]
\label{thm:TT}
Let $f$ be a function such that $\widetilde{f}^{(n+1)}(x)$
exists for all $x$ in a range that contains the number $a$.
Then, 
$$
f(x) = P_n(x) \oplus R_n(x)
$$
for all $x$, where
$$
P_n(x) = {\overset{n}{\widetilde{\sum_{k=0}}}} 
e^{\frac{1}{k!}} \odot \widetilde{f}^{(k)}(a)
\odot \left(x \ominus a\right)^{\{k\}}
$$
is the Taylor polynomial of degree $n$ and
$$
R_n(x) = e^{\frac{1}{(n+1)!}} \odot 
\widetilde{f}^{(n+1)}(c) 
\odot \left(x \ominus a\right)^{\{n+1\}}
$$
is the remainder term in Lagrange form,
for some number $c$ between $a$ and $x$.
\end{Theorem}

Suppose $f$ is a function that has derivatives
of all orders over an interval centered on $a$.
If $\lim_{n\rightarrow +\infty} R_n(x) = 1$
for all $x$ in the interval, then the Taylor
series is convergent and converges to $f(x)$:
$$
f(x) = {\overset{+\infty}{\widetilde{\sum_{k=0}}}} 
e^{\frac{1}{k!}} \odot \widetilde{f}^{(k)}(a) 
\odot \left(x \ominus a\right)^{\{k\}}.
$$

As examples of convergent series, one has:
\begin{equation*}
\begin{split}
e^x &= {\overset{+\infty}{\widetilde{\sum_{k=0}}}}
e^{\frac{1}{k!}} \odot x^{\{k\}},\\
\cos_e(x) &= {\overset{+\infty}{\widetilde{\sum_{k=0}}}}
e^{\frac{(-1)^k}{2k!}} \odot x^{\{2k\}}.
\end{split}
\end{equation*}

The Taylor's theorem given by Theorem~\ref{thm:TT}
has a natural extension for functions of several 
variables~\cite{MR3725737,e22101180}. Here we 
proceed by briefly reviewing integration.
For more on the alpha-arithmetic, its topology and analysis,
we refer the reader to the literature. For~example, 
mean value theorems can be found in the original book 
of Grossman and Katz of 1972~\cite{MR0430173}; 
for a recent reference with detailed proofs, see~\cite{MR3813861}.


\subsection{Integrals}

The notion of integral for the non-Newtonian calculus
under consideration is a type of product integration
~\cite{MR2917851}. As~expected, the~function $F(x)$ is
an antiderivative of function $f(x)$ on the interval $I$
if $\widetilde{F}(x) = f(x)$ for all $x \in I$. 
The~indefinite integral of $f(x)$ is denoted~by
$$
\strokedint f(x) \tilde{d}x = F(x) \oplus c,
$$
where $c$ is a constant. Examples are 
(see~\cite{MR695495,MR4188122,MR4247781}):
\begin{itemize}
\item If $k$ is a positive constant, then 
$\displaystyle \strokedint k \tilde{d}x = \left(k\odot x\right) 
\oplus c$---in particular, if~$k = 1$, 
then $\displaystyle \strokedint 1 \tilde{d}x = c$;
\item $\displaystyle \strokedint e^n \tilde{d}x = x^n \oplus c$;
\item $\displaystyle \strokedint x^{\{n\}} \tilde{d}x 
= e^{\frac{1}{n+1}} \odot x^{\{n+1\}} \oplus c$;
\item $\displaystyle \strokedint e \oslash x \, \tilde{d}x 
= \ln(x) \oplus c$;
\item $\displaystyle \strokedint e^{a\odot x} \tilde{d}x 
= \left(a^{\{-1\}}\odot e^{a\odot x}\right) \oplus c$;
\item $\displaystyle \strokedint e^{r x^2} \tilde{d}x 
= e^{\frac{r x^2}{2}} \oplus c$.
\end{itemize}

The definite integral of $f$ on $[a,b]$ is denoted by
$$
\strokedint_{a}^{b} f(x) \tilde{d}x.
$$

If $f$ is positive and continuous on $[a,b]$, then
$f$ is integrable in $[a,b]$. The~following properties~hold:
\begin{enumerate}
\item[(i)] $\displaystyle \strokedint_{a}^{b} f(x) \tilde{d}x 
= \strokedint_{a}^{c} f(x) \tilde{d}x \oplus
\strokedint_{c}^{b} f(x) \tilde{d}x$;
\item[(ii)] $\displaystyle \strokedint_{a}^{b} 
\left(f(x) \oplus g(x)\right) \tilde{d}x 
= \strokedint_{a}^{b} f(x) \tilde{d}x \oplus
\strokedint_{a}^{b} g(x) \tilde{d}x$;
\item[(iii)] $\displaystyle \strokedint_{a}^{b} 
\left(f(x) \ominus g(x)\right) \tilde{d}x 
= \strokedint_{a}^{b} f(x) \tilde{d}x \ominus
\strokedint_{a}^{b} g(x) \tilde{d}x$;
\item[(iv)] $\displaystyle \strokedint_{a}^{b} c \odot f(x) \tilde{d}x 
= c \odot \strokedint_{a}^{b} f(x) \tilde{d}x$;
\item[(v)] If $m \leq f(x) \leq M$ for all $x \in [a,b]$, then 
$m \odot \left(b \ominus a\right) \leq 
\displaystyle \strokedint_{a}^{b} f(x) \tilde{d}x 
\leq M \odot \left(b \ominus a\right)$;
\item[(vi)] If $f(x) \leq g(x)$ for all $x \in [a,b]$,
then $\displaystyle \strokedint_{a}^{b} f(x) \tilde{d}x 
\leq \strokedint_{a}^{b} g(x) \tilde{d}x$.
\end{enumerate}

If $f$ is positive and integrable on $[a,b]$, 
then $F$ defined on $[a,b]$ by
$$
F(x) = \strokedint_{a}^{x} f(t) \tilde{d}t 
$$
is continuous over $[a,b]$. Moreover, the~fundamental
theorems of integral calculus hold: if $f$ is continuous
in $x \in [a,b]$, then $F$ is differentiable at $x$ with
$$
\widetilde{F(x)} = f(x); 
$$
if $f = \widetilde{h}$ for some function $h$, then 
$$
\strokedint_{a}^{b} f(x) \tilde{d}x 
= h(b) \ominus h(a).
$$

For more on the $\alpha$-arithmetic, its generalized real analysis, its 
fundamental topological properties related to non-Newtonian metric 
spaces and its calculus, including non-Newtonian differential
equations and its applications, see~\cite{PAP2008368,MR3073475,%
duyar2015some,MR3463533,MR4104356,MR4188122,MR4170025}.
For gentle, thorough and modern introduction to the subject of
non-Newtonian calculi, we also refer the reader to the recent 
book~\cite{book2021:Burgin:Czachor}. 
Now we proceed with our original results.


\section{Results}
\label{sec:3}

In order to develop a non-Newtonian calculus of variations
(dynamic optimization), we begin by first proving some
necessary results of static~optimization.


\subsection{Static~Optimization}
\label{sec:3.1}

Given $\epsilon > 1$, let
 $$
\mathcal{B}(\bar{x},\epsilon) 
:= \left\{x \in \mathbb{R}^{+} : d(x,\bar{x}) \leq \epsilon\right\}
= \left\{x \in \mathbb{R}^{+} : [[x\ominus\bar{x}]] \leq \epsilon\right\}.
$$

Note that for $a,b \in \mathbb{R}^{+}$, one has
$$
a = b \Leftrightarrow \frac{a}{b} = 1 \ \left(\text{or }\frac{b}{a} = 1\right)
\Leftrightarrow a \ominus b = 1 \text{ (or $b \ominus a = 1$)}. 
$$

Similarly for inequalities, for~example,
$$
a < b \Leftrightarrow \frac{a}{b} < 1 \Leftrightarrow a \ominus b < 1.
$$

This means that
$$
\mathcal{B}(\bar{x},\epsilon) 
= \left\{x \in \mathbb{R}^{+} : d(x,\bar{x}) \ominus \epsilon \leq 1\right\}.
$$

\begin{Definition}[local minimizer]
Let $f : (a,b) \rightarrow \mathbb{R}^{+}$ and consider the problem
of minimizing $f(x)$, $x \in (a,b)$. We say that $x \in (a,b)$ 
is a (local) minimizer of $f$ in $(a,b)$ if there exists $\epsilon > 1$
such that $f(x) \leq f(y)$ (i.e., $f(x) \ominus f(y) \leq 1$) for all
$y \in \mathcal{B}(x, \epsilon) \cap (a,b)$. In~this case, we say that
$f(x)$ is a (local) minimum.
\end{Definition}

Another important concept in optimization is that of descent~direction.

\begin{Definition}[descent direction]
A $d \in \mathbb{R}^{+}$ is said to be a descent direction of $f$ at $x$
if $f\left(x\oplus \epsilon \odot d\right) < f(x)$ $\forall$ $\epsilon > 1$
sufficiently close to 1, or~equivalently, if~$f\left(x\oplus \epsilon 
\odot d\right) \ominus f(x) < 1$ for all $\epsilon > 1$
sufficiently close to 1.
\end{Definition}

\begin{Remark}
From the chain rule and other properties of Section~\ref{sec:2},
it follows that
$$
\frac{\tilde{d}}{\tilde{d}x}\left[f\left(x\oplus \epsilon \odot d\right)\right]
= \widetilde{f}\left(x\oplus \epsilon \odot d\right)
$$
and
\begin{equation}
\label{eq:my:1}
\frac{\tilde{d}}{\tilde{d}\epsilon}\left[f\left(x\oplus \epsilon \odot d\right)\right]
= \widetilde{f}\left(x\oplus \epsilon \odot d\right)\odot d.
\end{equation}

In particular, we get from \eqref{eq:my:1} that
\begin{equation*}
\left.\frac{\tilde{d}}{\tilde{d}\epsilon}\left[f\left(x\oplus 
\epsilon \odot d\right)\right]\right|_{\epsilon = 1}
= \widetilde{f}(x)\odot d.
\end{equation*}
\end{Remark}

Our first result allow us to identify a descent direction of $f$
at $x$ based on the derivative of $f$ at $x$.

\begin{Theorem}
\label{thm:1}
Let $f$ be differentiable. If~there exists $d \in \mathbb{R}^{+}$
such that $\widetilde{f}(x)\odot d < 1$, then $d$ is a descent direction
of $f$ at $x$.
\end{Theorem}

\begin{proof}
We know from Taylor's theorem (Theorem~\ref{thm:TT}) that
\begin{equation}
\label{eq:my:2}
f\left(x\oplus \epsilon \odot d\right) 
= f(x) \oplus \epsilon \odot \widetilde{f}(x) \odot d
\oplus R_1\left(x\oplus \epsilon \odot d\right),
\end{equation}
where 
\begin{equation*}
\begin{split}
R_1\left(x\oplus \epsilon \odot d\right)
&= \left(\widetilde{f}^{(2)}(c)\right)^{\frac{1}{2}}
\odot \left(\epsilon \odot d\right)^{\{2\}}\\
&= \epsilon^{\{2\}} \odot \left(\widetilde{f}^{(2)}(c)\right)^{\frac{1}{2}}
\odot d^{\{2\}}
\end{split}
\end{equation*}
with $c$ being in the interval between $x$ and $x\oplus \epsilon \odot d$.
The equality \eqref{eq:my:2} can be written in the following equivalent form:

\begin{equation}
\label{eq:my:3}
f\left(x\oplus \epsilon \odot d\right) 
\ominus f(x) = \epsilon \odot \widetilde{f}(x) \odot d
\oplus R_1\left(x\oplus \epsilon \odot d\right).
\end{equation}

Recalling that $a\odot b^{\{-1\}} = a \oslash b$,
$a\odot a^{\{-1\}} = e$ and $\odot$ is distributive
over $\oplus$, we get from \eqref{eq:my:3} that
\begin{equation}
\label{eq:my:4}
\begin{split}
\left(f\left(x\oplus \epsilon \odot d\right) \ominus f(x)\right)\oslash \epsilon 
&= \widetilde{f}(x) \odot d \oplus \epsilon^{\{-1\}}\odot R_1\left(x\oplus \epsilon \odot d\right)\\
&= \widetilde{f}(x) \odot d \oplus \epsilon \odot \left(\widetilde{f}^{(2)}(c)\right)^{\frac{1}{2}}
\odot d^{\{2\}}.
\end{split}
\end{equation}

Now we note that as $\epsilon\rightarrow 1$, one has
$\epsilon \odot \left(\widetilde{f}^{(2)}(c)\right)^{\frac{1}{2}}\odot d^{\{2\}}\rightarrow 1$
so that the right-hand side of \eqref{eq:my:4} converges to $\widetilde{f}(x) \odot d$.
From the hypothesis $\widetilde{f}(x)\odot d < 1$ of our theorem, this means that,
for $\epsilon > 1$ sufficiently close to 1, the~right-hand side of \eqref{eq:my:4}
is strictly less than one. Thus, for~$\epsilon > 1$ sufficiently close to 1,
\begin{equation}
\label{eq:my:5}
\left(f\left(x\oplus \epsilon \odot d\right) \ominus f(x)\right)\oslash \epsilon < 1.
\end{equation}

Recalling that $a \oslash \epsilon < 1 \Leftrightarrow a^{1/\ln(\epsilon)} < 1$,
we conclude from \eqref{eq:my:5} that for $\epsilon$ sufficiently close to 1
we have $f\left(x\oplus \epsilon \odot d\right) \ominus f(x) < 1$;
that is, $d$ is a descent direction of $f$ at $x$.
\end{proof}

As a corollary of Theorem~\ref{thm:1}, we obtain 
Fermat's necessary optimality condition, which
gives us a method to find local minimizers (or maximizers)
of differentiable functions on open sets, by~showing that 
every local extremizer of the function is a stationary point 
(the non-Newtonian function's derivative is one at that point).

\begin{Theorem}[Fermat's theorem--stationary points]
\label{thm:2}
Let $f : (a,b) \rightarrow \mathbb{R}^{+}$
be differentiable. If~$x \in (a,b)$ is a minimizer
of $f$, then $\widetilde{f}(x) = 1$.
\end{Theorem}

\begin{proof}
We want to prove that for a minimizer $x$ we must have
\mbox{$\widetilde{f}(x) = 1 \Leftrightarrow 1 \ominus \widetilde{f}(x) = 1$}.
We do the proof by contradiction. Assume that $\widetilde{f}(x) \neq 1$;
that is, $1 \ominus \widetilde{f}(x) \neq 1$.
Let \mbox{$d = 1 \ominus \widetilde{f}(x)$}.~Then,
\begin{equation*}
\begin{split}
\widetilde{f}(x) \odot d 
&= \widetilde{f}(x) \odot \left(1 \ominus \widetilde{f}(x)\right)
= 1 \ominus \widetilde{f}(x)^{\{2\}}\\
&= \frac{1}{\widetilde{f}(x) \odot \widetilde{f}(x)}
= \left(\frac{1}{\widetilde{f}(x)}\right)^{\ln\left(\widetilde{f}(x)\right)}
\end{split}
\end{equation*}
and since $g(y) = \left(\frac{1}{y}\right)^{\ln(y)}$ is a function with
$0 < g(y) < 1$ for all $y \neq 1$, we conclude that $\widetilde{f}(x) \odot d < 1$.
It follows from Theorem~\ref{thm:1} that $d$ is a descent direction of $f$ at $x$,
and therefore, from~the definition of descent direction, $x$ is not a local minimizer.
\end{proof}

In the next section, we make use of Theorem~\ref{thm:2}
to prove the non-Newtonian Euler--Lagrange~equation.


\subsection{Dynamic~Optimization}
\label{sec:3.2}

A central tool in dynamic optimization, both in the
calculus of variations and optimal control~\cite{MR641031},
is integration by parts. In~what follows, we use the following
notation:
$$
\psi(x)\big|_{a}^{b} = \psi(b) \ominus \psi(a).
$$

\begin{Theorem}[integration by parts]
\label{thm:3}
Let $f : [a,b] \rightarrow \mathbb{R}^{+}$
and $g : [a,b] \rightarrow \mathbb{R}^{+}$
be differentiable. The~following formula
of integration by parts holds:

\begin{equation}
\label{eq:my:6}
\strokedint_{a}^{b} \widetilde{f}(x) \odot g(x) \tilde{d}x
= f(x) \odot g(x)\big|_{a}^{b}  \ominus
\strokedint_{a}^{b} f(x) \odot \widetilde{g}(x) \tilde{d}x.
\end{equation}
\end{Theorem}

\begin{proof}
From the derivative of a product, we know that
\begin{equation}
\label{eq:my:7}
\frac{\tilde{d}}{\tilde{d}x} \left[f(x) \odot g(x)\right]
= \widetilde{f}(x) \odot g(x) \oplus f(x) \odot \widetilde{g}(x).
\end{equation}

On the other hand, the~fundamental theorem of integral calculus tell us that
$$
\strokedint_{a}^{b} \frac{\tilde{d}}{\tilde{d}x} 
\left[f(x) \odot g(x)\right] \tilde{d}x
= f(x) \odot g(x)\big|_{a}^{b}.
$$

Therefore, by~integrating \eqref{eq:my:7} from $a$ to $b$, we conclude that
$$
\strokedint_{a}^{b} \left[\widetilde{f}(x) \odot g(x) \right]\tilde{d}x
\oplus 
\strokedint_{a}^{b} \left[f(x) \odot \widetilde{g}(x)\right]\tilde{d}x
= f(x) \odot g(x)\big|_{a}^{b}\, ,
$$
which is equivalent to \eqref{eq:my:6}.
\end{proof}

We are now in a condition to formulate the fundamental problem of the
calculus of variations: to minimize the integral functional
$$
\mathscr{F}[y] = \strokedint_{a}^{b} 
L\left(x,y(x),\widetilde{y}(x)\right) \tilde{d}x
$$
over all smooth functions $y$ on $[a,b]$
with fixed end points $y(a) = y_a$ and $y(b) = y_b$.  
The central result of the calculus of variations and classical mechanics 
is the celebrated Euler--Lagrange equation, 
whose solutions are stationary points of the given action functional.
We restrict ourselves here to the classical framework 
of the calculus of variations, where both the Lagrangian 
$L$ and admissible functions $y$ are  smooth enough: 
typically, one considers $L \in C^2$  and $y \in C^2$, so that
one can look to the Euler--Lagrange equation as
a second-order ordinary differential equation~\cite{MR2004181}.
We adopt such assumptions here. We denote the problem by $(P)$.

Before proving the Euler--Lagrange equation (the necessary optimality
condition for problem $(P)$), we first need to prove a non-Newtonian
analogue of the fundamental lemma of the calculus of~variations.

\begin{Lemma}[fundamental lemma of the calculus of variations]
\label{FL}	
If $f : [a,b] \rightarrow \mathbb{R}^{+}$ is a positive continuous
function such that
\begin{equation}
\label{eq:my:8}
\strokedint_{a}^{b} f(x) \odot h(x) \tilde{d}x = 1
\end{equation}
for all functions $h(x)$ that are continuous for $a\leq x \leq b$
with $h(a) = h(b) = 1$, then $f(x) = 1$ for all $x \in [a,b]$.
\end{Lemma}

\begin{proof}
We do the proof by contradiction.
Suppose the function $f$ is not one---$f(x) > 1$---at some point 
$x \in [a,b]$. Then, by~continuity, $f(x) > 1$ for all $x$ 
in some interval $[x_1, x_2] \subset [a,b]$. If we set
$$
h(x) =
\begin{cases}
\left(x \ominus x_1\right) \odot \left(x_2 \ominus x\right)
& \text{if } x \in [x_1,x_2],\\
1 & \text{if } x \in [a,b] \setminus [x_1,x_2],
\end{cases}
$$
then $h(x)$ satisfies the assumptions of the lemma; i.e.,
$h(x)$ is continuous for $x \in [a,b]$ with
$h(a) = 1$ and $h(b) = 1$. We have

\begin{equation}
\label{eq:my:9}
\begin{split}
\strokedint_{a}^{b} f(x) \odot h(x) \tilde{d}x 
&= \strokedint_{a}^{x_1} 1 \tilde{d}x 
\oplus \strokedint_{x_1}^{x_2} f(x) 
\odot \left(x \ominus x_1\right) 
\odot \left(x_2 \ominus x\right) \tilde{d}x 
\oplus \strokedint_{x_2}^{b} 1 \tilde{d}x \\
&= \strokedint_{x_1}^{x_2} f(x) 
\odot \left(x \ominus x_1\right) 
\odot \left(x_2 \ominus x\right) \tilde{d}x. 
\end{split}
\end{equation}

Let us analyze the integrand $\beta(x)$ of \eqref{eq:my:9}:
\begin{equation*}
\begin{split}
\beta(x) 
&= f(x) \odot \left[
\left(x \ominus x_1\right) 
\odot \left(x_2 \ominus x\right)\right]
= f(x) \odot \left[ \frac{x}{x_1} \odot \frac{x_2}{x}\right]\\
&= f(x) \odot \left[\left(\frac{x}{x_1}\right)^{\ln\left(\frac{x_2}{x}\right)}\right].
\end{split}
\end{equation*}

Since $f(x) > 1$ and 
$\left(\frac{x}{x_1}\right)^{\ln\left(\frac{x_2}{x}\right)} > 1$
for any $x \in (x_1,x_2)$, we also have that $\beta(x) > 1$
for any $x \in (x_1,x_2)$, with~$\beta(x) = 1$ at $x = x_1$
and $x = x_2$. It follows that
$$
\strokedint_{a}^{b} f(x) \odot h(x) \tilde{d}x 
= \strokedint_{x_1}^{x_2} \beta(x) \tilde{d}x 
> 1 \odot \left(x_2 \ominus x_1\right) = 1.
$$

This contradicts \eqref{eq:my:8} and proves the lemma.
\end{proof}

Now we formulate and prove the analog of the Euler--Lagrange
differential equation for our problem $(P)$.

\begin{Theorem}[Euler--Lagrange equation]
\label{thm:4}	
If $y(x)$, $x \in [a,b]$, is a solution to the problem
\begin{equation}
\label{P}
\tag{$P$}
\mathscr{F}[y] = \strokedint_{a}^{b}
L\left(x,y(x),\widetilde{y}(x)\right) \tilde{d}x 
\longrightarrow \min_{y \in \mathcal{Y}\left(y_a;y_b\right)}
\end{equation}
with
\begin{equation*}
\mathcal{Y}\left(y_a;y_b\right) :=
\left\{y \in C^2\left([a,b];\mathbb{R}^+\right) : 
y(a) = y_a,\  y(b) = y_b,\  y(x) > 0\  \forall \ x \in [a,b] \right\},  
\end{equation*}
then $y(x)$ satisfies the Euler--Lagrange equation
\begin{equation}
\label{ELeq}
\widetilde{L}_{y}\left(x,y(x),\widetilde{y}(x)\right)
= \frac{\tilde{d}}{\tilde{d}x} 
\widetilde{L}_{\widetilde{y}}\left(x,y(x),\widetilde{y}(x)\right) 
\end{equation}
for all $x \in [a,b]$.
\end{Theorem}

\begin{proof}
Let $y(x)$, $x \in [a,b]$, be a minimizer of \eqref{P}.
Then, function $\left(y \oplus \epsilon \odot h\right)(x)$,
$x\in [a,b]$, belongs to $\mathcal{Y}\left(y_a;y_b\right)$ 
for any function $h \in \mathcal{Y}\left(1;1\right)$
and for any $\epsilon$ in an open neighborhood of $1$. Note that
$\left(y \oplus \epsilon \odot h\right)(x) = y(x)$
for $\epsilon = 1$. This means that for any smooth function 
$h(x)$, $x \in [a,b]$, satisfying $h(a) = h(b) =1$, 
the function $\varphi(\epsilon)$ defined by
\begin{equation}
\label{eq:my:11}
\begin{split}
\varphi(\epsilon)
&= \mathscr{F}[y \oplus \epsilon \odot h] 
= \strokedint_{a}^{b}
L\left(x,\left(y \oplus \epsilon \odot h\right)(x),
\widetilde{\left(y \oplus \epsilon \odot h\right)}(x)\right) \tilde{d}x\\
&= \strokedint_{a}^{b}
L\left(x,y(x) \oplus \epsilon \odot h(x),
\widetilde{y}(x) \oplus \epsilon \odot \widetilde{h}(x)\right) \tilde{d}x
\end{split}
\end{equation}
has a minimizer for $\epsilon = 1$. It follows from Fermat's theorem
(Theorem~\ref{thm:2}) that
$$
\widetilde{\varphi}(1) = 
\left.\frac{\tilde{d}}{\tilde{d}\epsilon}\varphi(\epsilon)\right|_{\epsilon=1} = 1.
$$

By differentiating \eqref{eq:my:11} with respect to $\epsilon$, 
and then putting $\epsilon = 1$, we get from the chain rule 
and relations
$$
\frac{\tilde{d}}{\tilde{d}\epsilon} 
\left(y(x) \oplus \epsilon \odot h(x)\right) 
= h(x), \quad 
\frac{\tilde{d}}{\tilde{d}\epsilon} 
\left(\widetilde{y}(x) \oplus \epsilon \odot \widetilde{h}(x)\right) 
= \widetilde{h}(x), 
$$
that
\begin{equation}
\label{eq:my:12}
1 = \strokedint_{a}^{b} \left[ 
\widetilde{L}_{y}\left(x,y(x),\widetilde{y}(x)\right) \odot h(x) 
\oplus \widetilde{L}_{\widetilde{y}}\left(x,y(x),\widetilde{y}(x)\right) 
\odot \widetilde{h}(x)\right] \tilde{d}x.
\end{equation}

From integration by parts (Theorem~\ref{thm:3}),
and the fact that $h(a) = 1$ and $h(b) = 1$, one has
\begin{multline*}
\strokedint_{a}^{b}
\widetilde{L}_{\widetilde{y}}\left(x,y(x),\widetilde{y}(x)\right) 
\odot \widetilde{h}(x) \tilde{d}x
= \left. h(x) \odot \widetilde{L}_{\widetilde{y}}\left(x,y(x),
\widetilde{y}(x)\right)\right|_{a}^{b}\\
\ominus \strokedint_{a}^{b}
\frac{\tilde{d}}{\tilde{d}x} \left[
\widetilde{L}_{\widetilde{y}}\left(x,y(x),\widetilde{y}(x)\right)\right] 
\odot h(x)\tilde{d}x,
\end{multline*}
that is,
\begin{equation}
\label{eq:my:13}
\strokedint_{a}^{b}
\widetilde{L}_{\widetilde{y}}\left(x,y(x),\widetilde{y}(x)\right) 
\odot \widetilde{h}(x) \tilde{d}x
= 1 \ominus \strokedint_{a}^{b}
\frac{\tilde{d}}{\tilde{d}x} \left[
\widetilde{L}_{\widetilde{y}}\left(x,y(x),\widetilde{y}(x)\right)\right] 
\odot h(x)\tilde{d}x.
\end{equation}

Using equality \eqref{eq:my:13} in the necessary condition \eqref{eq:my:12},
we get that
\begin{multline}
\label{eq:my:14}
1 = \strokedint_{a}^{b}
\widetilde{L}_{y}\left(x,y(x),\widetilde{y}(x)\right) \odot h(x) 
\tilde{d}x
\oplus 
1 \ominus \strokedint_{a}^{b}
\frac{\tilde{d}}{\tilde{d}x} \left[
\widetilde{L}_{\widetilde{y}}\left(x,y(x),\widetilde{y}(x)\right)\right] 
\odot h(x)\tilde{d}x\\
\Leftrightarrow
\strokedint_{a}^{b} \left[
\widetilde{L}_{y}\left(x,y(x),\widetilde{y}(x)\right) 
\ominus \frac{\tilde{d}}{\tilde{d}x} \left(
\widetilde{L}_{\widetilde{y}}\left(x,y(x),\widetilde{y}(x)\right)\right)
\right] \odot h(x) \tilde{d}x = 1.
\end{multline}

The result follows from the fundamental lemma of the calculus of variations
(\mbox{Lemma~\ref{FL}}) applied to \eqref{eq:my:14}:
$\widetilde{L}_{y}\left(x,y(x),\widetilde{y}(x)\right) 
\ominus \frac{\tilde{d}}{\tilde{d}x} \left(
\widetilde{L}_{\widetilde{y}}\left(x,y(x),
\widetilde{y}(x)\right)\right) = 1$ for all $x \in [a,b]$.
\end{proof}

To illustrate our main result, let us see an example.
Consider the following problem of the calculus of variations:
\begin{equation}
\label{eq:my:18}
\begin{gathered}
\mathscr{F}[y] = \sqrt{e} \odot 
\strokedint_{1}^{e^{2\pi}}
\left[ \widetilde{y}^{\{2\}}(x) \ominus y^{\{2\}}(x)\right] 
\tilde{d}x \longrightarrow \min,\\
y(1) = e, \quad y(e^{2\pi}) = e^{-1}.
\end{gathered}
\end{equation}

Theorem~\ref{thm:4} tell us that the solution of \eqref{eq:my:18}
must satisfy the Euler--\mbox{Lagrange Equation~\eqref{ELeq}}. In~this example,
the Lagrangian $L$ is given by $L(x,y,\widetilde{y}) 
= \sqrt{e} \odot \left(\widetilde{y}^{\{2\}} \ominus y^{\{2\}}\right)$,
so that
\begin{equation}
\label{eq:my:19}
\begin{split}
\widetilde{L}_{y}\left(x,y,\widetilde{y}\right)
&= \sqrt{e} \odot \left(1 \ominus e^2 \odot y\right),\\
\widetilde{L}_{\widetilde{y}}\left(x,y,\widetilde{y}\right)
&=\sqrt{e} \odot \left(e^2 \odot \widetilde{y} \ominus 1\right)
= \sqrt{e} \odot \left(e^2 \odot \widetilde{y}\right).
\end{split}
\end{equation}

Noting that $\sqrt{e} = e \oslash e^2 = \left(e^2\right)^{\{-1\}}$, 
the~equalities \eqref{eq:my:19} simplify to
\begin{equation*}
\widetilde{L}_{y}\left(x,y,\widetilde{y}\right)
= 1 \ominus y,\quad
\widetilde{L}_{\widetilde{y}}\left(x,y,\widetilde{y}\right)
= \widetilde{y}
\end{equation*}
and the Euler--Lagrange Equation \eqref{ELeq} takes the form
\begin{equation}
\label{eq:my:10}
1 \ominus y(x) = \widetilde{y}^{(2)}(x)
\Leftrightarrow \widetilde{y}^{(2)}(x) \oplus y(x) = 1.
\end{equation}

The second-order differential Equation \eqref{eq:my:10}
has solutions of the form
$$
y(x) = c_1 \odot \cos_e(x) \oplus c_2 \odot \sin_e(x),
$$
where $c_1$ and $c_2$ are constants. Given the boundary
conditions $y(1) = e$ and $y\left(e^{2\pi}\right) = e^{-1}$,
we conclude that the Euler--Lagrange extremal for problem
\eqref{eq:my:18} is given by
$$
y(x) = e \odot \cos_e(x) \oplus e^2 \odot \sin_e(x)
\Leftrightarrow
y(x) = e^{2\sin\left(\ln(x)\right) + \cos\left(\ln(x)\right)}.
$$


\section{Discussion}
\label{sec:4}

One can say that the calculus of variations began 
in 1687 with Newton's minimal resistance 
problem~\cite{MR2164553,MR2323273,MR2201075}. It immediately occupied 
the attention of Bernoulli (\mbox{1655--1705}), but~it was Euler (1707--1783) 
who first elaborated mathematically on the subject, beginning in 1733. 
Lagrange (1736--1813) was influenced by Euler's work and contributed 
significantly to the theory, introducing a purely analytic approach 
to the subject based on additive variations $y + \epsilon h$, 
whose essence we still follow today~\cite{MR2004181}. The~calculus 
of variations is concerned with the minimization of integral functionals:
\begin{equation}
\label{eq:D1}
\int_{a}^{b} L\left(x,y(x),y'(x)\right) dx \longrightarrow \min.
\end{equation}

However, as~observed in~\cite{MR2378083}, there are other
interesting problems that arise in applications in which
the functionals to be minimized are not of the form \eqref{eq:D1}.
For example, the~planning of a firm trying to program its
production and investment policies to reach a given production
rate and to maximize its future market competitiveness at a
given time horizon, can be mathematically stated in the form
\begin{equation}
\label{eq:D2}
\left(\int_{a}^{b} L_1\left(x,y(x),y'(x)\right) dx\right) 
\cdot
\left(\int_{a}^{b} L_2\left(x,y(x),y'(x)\right) dx\right) 
\longrightarrow \min
\end{equation}
\textls[-15]{(see~\cite{MR2378083}). 
Another example, also given in~\cite{MR2378083},
appears when dealing with the so called ``slope stability problem'',
which is described mathematically as minimizing a quotient~functional:}
\begin{equation}
\label{eq:D3}
\frac{\displaystyle \int_{a}^{b} L_1\left(x,y(x),y'(x)\right) 
dx}{\displaystyle \int_{a}^{b} L_2\left(x,y(x),y'(x)\right) dx} 
\longrightarrow \min.
\end{equation}

Such multiplicative integral minimization problems 
that arise in different applications are nonstandard 
problems of the calculus of variations, but~they can
be naturally modeled in the non-Newtonian calculus
of variations, and~then solved in a rather standard
way, using non-Newtonian Lagrange variations of the 
form $y \oplus \epsilon \odot h$,
as we have proposed here. Therefore, we claim that 
the non-Newtonian calculus of variations
just introduced may be useful for dealing with
multiplicative functionals that arise
in economics, physics and~biology.

In this paper we have restricted ourselves to the ideas
and to the central results of any calculus of variations:
the celebrated Euler--Lagrange equation, which is a
first-order necessary optimality condition. Of~course
our results can be extended, for~example, by~relaxing
the considered hypotheses and enlarging the space
of admissible functions, which we have considered here 
to be $C^2$, or~considering vector functions instead 
of scalar ones. We leave such generalizations
to the interested and curious reader. In~fact,
much remains now to be done. As~possible future 
research directions we can mention:  obtaining
natural boundary conditions (sometimes also called
transversality conditions) to be satisfied at a boundary
point $a$ and/or $b$, when $y(a)$ and/or $y(b)$ are free
or restricted to take values on a given curve; obtaining
second-order necessary conditions; obtaining sufficient conditions;
to investigate nonadditive isoperimetric problems; etc.


\section{Conclusions}
\label{sec:5}

\textls[-5]{In this work, a~new calculus of variations was proposed, based on the non-Newtonian 
approach introduced by Grossman and Katz, thereby avoiding problems about nonnegativity. 
A new relation was proved, the~multiplicative Euler--Lagrange  differential 
Equation \eqref{ELeq}, which each solution of a non-Newtonian variational problem, 
with admissible functions taking positive values only, must satisfy.
An example was provided for illustration~purposes.}

Grossman and Katz have shown that infinitely many 
calculi can be constructed independently~\cite{MR0430173}.
Each of these calculi provide different perspectives 
for approaching many problems in science 
and engineering~\cite{MR2588234}. Additionally, 
a mathematical problem, which is difficult or 
impossible to solve in one calculus, can be easily 
revealed through another calculus~\cite{MR3794497,MR4247781}.

Since the pioneering work of Grossman and Katz, non-Newtonian
calculi have been a topic for new study areas in mathematics
and its applications~\cite{MR3461691,MR3528851}. Particularly,
Stanley~\cite{stanley1999multiplicative},
C{\'o}rdova-Lepe~\cite{cordova2006multiplicative},
Slav\'{\i}k~\cite{MR2917851}, Pap~\cite{PAP2008368}
and Bashirov~et~al.~\cite{MR2356052}, have called
the attention of mathematicians to the topic. 
More recently, non-Newtonian calculi have become 
a hot topic in economic and finance~\cite{MR3176131}, 
quantum calculus~\cite{MR3461691},
complex analysis~\cite{MR2734313,MR3198166,MR3891225}, 
numerical analysis~\cite{MR3725737,MR4164307,MR4222466},
inequalities~\cite{MR3016664,MR4016182}, 
biomathematics~\cite{MR2878647,MR4188122}
and mathematical education~\cite{MR4247781}.

Here we adopted the non-Newtonian calculus as originally introduced
by Grossman and Katz~\cite{MR0430173,MR557734,MR695495}
and recently developed by C{\'o}rdova-Lepe~\cite{cordova2006multiplicative,MR2724186}
and collaborators~\cite{MORA20121245}: see the recent reviews
in~\cite{MR4188122,MR4247781}. Roughly speaking, the~key to 
understanding such calculus, valid for positive functions,
is a formal substitution, where one replaces addition 
and subtraction with multiplication and division, respectively;  
multiplication in standard calculus is replaced by exponentiation 
in the non-Newtonian case, and~thus, division by exponentiation 
with the reciprocal exponent. Our main contribution here was to develop, 
for the first time in the literature, 
a suitable non-Newtonian calculus of variations
that  minimizes a non-Newtonian integral functional
with a Lagrangian that depends on the non-Newtonian derivative.
The main result is a first-order necessary optimality condition
of Euler--Lagrange~type. 

We trust that the present paper marks the beginning of a fruitful 
road for Non-Newtonian (NN) mechanics, NN calculus of variations
and NN optimal control, thereby calling attention to and serving 
as inspiration for a new generation of researchers. Currently, 
we are investigating the validity of Emmy Noether's principle 
in the NN/multiplicative calculus of variations here~introduced.


\vspace{6pt} 


\funding{This research was funded by 
The Center for Research and Development 
in Mathematics and Applications (CIDMA)
through the Portuguese Foundation for Science and Technology 
(FCT--Funda\c{c}\~{a}o para a Ci\^{e}ncia 
e a Tecnologia), grant number UIDB/04106/2020.}

\institutionalreview{Not applicable.}

\informedconsent{Not applicable.}

\dataavailability{Not applicable.}

\acknowledgments{The author is grateful
to Luna Shen, Managing Editor of Axioms,
for proposing a volume/special issue dedicated 
to his 50th anniversary, and~to Nat\'{a}lia Martins,
Ricardo Almeida, Cristiana J. Silva and
M. Rchid Sidi Ammi, who kindly accepted 
the invitation of Axioms to lead the~project.}

\conflictsofinterest{The author declares no conflict of interest. 
The funders had no role in the design of the study; 
in the collection, analyses, or~interpretation of data; 
in the writing of the manuscript
or in the decision to publish the~results.} 


\end{paracol}

\reftitle{References}


\section*{Short Biography of Author}

\bio
{\raisebox{0.95cm}{\includegraphics[width=3.5cm,height=5.3cm,clip,keepaspectratio]{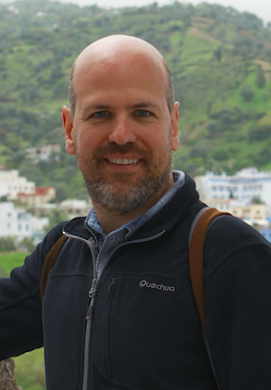}}}
{\textbf{Delfim Fernando Marado Torres} is a Portuguese Mathematician 
born 16 August 1971 in Nampula, Portuguese Mozambique. 
He obtained a PhD in Mathematics from the University of Aveiro (UA) in 2002, 
and habilitation in Mathematics, UA, in 2011. He is a full professor 
of mathematics since \mbox{9 March 2015}. He has been the Director of the R\&D Unit CIDMA, 
the largest Portuguese research center for mathematics, 
and Coordinator of its Systems and Control Group. His main research areas
are calculus of variations and optimal control; optimization; fractional derivatives 
and integrals; dynamic equations on time scales; and mathematical biology. 
Torres has written outstanding scientific and pedagogical publications. In particular, 
he has co-authored two books with Imperial College Press and three books with Springer. 
He has strong experience in graduate and post-graduate student supervision 
and teaching in mathematics. Twenty PhD students in mathematics 
have successfully finished under his supervision. Moreover, he has been the leading member 
in several national and international R\&D projects, including EU projects 
and networks. Professor Torres has been, since 2013, the Director of the Doctoral Programme 
Consortium in Mathematics and Applications (MAP-PDMA) of Universities 
of Minho, Aveiro, and Porto. Delfim married 2003 and has one daughter and two~sons.}



\begin{thebibliography}{999}

\bibitem[Grossman and Katz(1972)]{MR0430173}
Grossman, M.; Katz, R.
\newblock {\em Non-{N}ewtonian Calculus}; Lee Press: Pigeon Cove, MA, USA,  1972.

\bibitem[Boruah \em{et~al.}(2021)Boruah, Hazarika, and Bashirov]{MR4164307}
Boruah, K.; Hazarika, B.; Bashirov, A.E.
\newblock Solvability of bigeometric differential equations by numerical
methods.
\newblock {\em Bol. Soc. Parana. Mat.} {\bf 2021}, {\em 39},~203--222.

\bibitem[Bashirov \em{et~al.}(2011)Bashirov, M\i~s\i rl\i, Tando\u{g}du, and
\"{O}zyap\i c\i]{MR2864779}
Bashirov, A.E.; M\i~s\i rl\i, E.; Tando\u{g}du, Y.; \"{O}zyap\i c\i, A.
\newblock On modeling with multiplicative differential equations.
\newblock {\em Appl. Math. J. Chin. Univ. Ser. B} {\bf 2011}, 
{\em 26},~425--438, doi:10.1007/s11766-011-2767-6.

\bibitem[Mora \em{et~al.}(2012)Mora, C\'{o}rdova-Lepe, and
Del-Valle]{MORA20121245}
Mora, M.; C\'{o}rdova-Lepe, F.; Del-Valle, R.
\newblock A non-{N}ewtonian gradient for contour detection in images with
multiplicative noise.
\newblock {\em Pattern Recognit. Lett.} {\bf 2012}, {\em 33},~1245--1256, 
doi:10.1016/j.patrec.2012.02.012.

\bibitem[Ozyapici and Bilgehan(2016)]{MR3452941}
Ozyapici, A.; Bilgehan, B.
\newblock Finite product representation via multiplicative calculus and its
applications to exponential signal processing.
\newblock {\em Numer. Algorithms} {\bf 2016}, {\em 71},~475--489, 
doi:10.1007/s11075-015-0004-8.

\bibitem[Czachor(2020)]{e22101180}
Czachor, M.
\newblock Unifying Aspects of Generalized Calculus.
\newblock {\em Entropy} {\bf 2020}, {\em 22}, 1180, doi:10.3390/e22101180.

\bibitem[Pinto \em{et~al.}(2020)Pinto, Torres, Campillay-Llanos, and
Guevara-Morales]{MR4188122}
Pinto, M.; Torres, R.; Campillay-Llanos, W.; Guevara-Morales, F.
\newblock Applications of proportional calculus and a non-{N}ewtonian logistic
growth model.
\newblock {\em Proyecciones} {\bf 2020}, {\em 39},~1471--1513.

\bibitem[\"{O}zyap\i c\i \em{et~al.}(2014)\"{O}zyap\i c\i, Riza, Bilgehan, and
Bashirov]{MR3271405}
\"{O}zyap\i c\i, A.; Riza, M.; Bilgehan, B.; Bashirov, A.E.
\newblock On multiplicative and {V}olterra minimization methods.
\newblock {\em Numer. Algorithms} {\bf 2014}, {\em 67},~623--636, 
doi:10.1007/s11075-013-9813-9.

\bibitem[Ja\"{e}ck(2019)]{MR3992455}
Ja\"{e}ck, F.
\newblock Calcul diff\'{e}rentiel et int\'{e}gral adapt\'{e} 
aux substitutions par {V}olterra.
\newblock {\em Hist. Math.} {\bf 2019}, {\em 48},~29--68, 
doi:10.1016/j.hm.2018.12.001.

\bibitem[Slav\'{\i}k(2007)]{MR2917851}
Slav\'{\i}k, A.
\newblock {\em Product Integration, Its History and Applications}; 
Matfyzpress: Prague, Czech Republic, 2007; p. iv+147.

\bibitem[C{\'o}rdova-Lepe(2006)]{cordova2006multiplicative}
C{\'o}rdova-Lepe, F.
\newblock The multiplicative derivative as a measure of elasticity in economics.
\newblock {\em TMAT Rev. Latinoam. Cienc. E Ing.} {\bf 2006}, {\em 2},~1--8.

\bibitem[Bashirov \em{et~al.}(2008)Bashirov, Kurp\i~nar, and \"{O}zyap\i c\i]{MR2356052}
Bashirov, A.E.; Kurp\i~nar, E.M.; \"{O}zyap\i c\i, A.
\newblock Multiplicative calculus and its applications.
\newblock {\em J. Math. Anal. Appl.} {\bf 2008}, {\em 337},~36--48, 
doi:10.1016/j.jmaa.2007.03.081.

\bibitem[Boruah and Hazarika(2017)]{MR3722924}
Boruah, K.; Hazarika, B.
\newblock On some generalized geometric difference sequence spaces.
\newblock {\em Proyecciones} {\bf 2017}, {\em 36},~373--395, 
doi:10.4067/s0716-09172017000300373.

\bibitem[Campillay-Llanos \em{et~al.}(2021)Campillay-Llanos, Guevara, Pinto,
and Torres]{MR4247781}
Campillay-Llanos, W.; Guevara, F.; Pinto, M.; Torres, R.
\newblock Differential and integral proportional calculus: How to find a
primitive for {$f(x)=1/sqrt{2\pi}e^{-(1/2)x^2}$}.
\newblock {\em Int. J. Math. Ed. Sci. Technol.} {\bf 2021}, 
{\em 52},~463--476, doi:10.1080/0020739X.2020.1763489.

\bibitem[Gurefe \em{et~al.}(2016)Gurefe, Kadak, Misirli, and Kurdi]{MR3528851}
Gurefe, Y.; Kadak, U.; Misirli, E.; Kurdi, A.
\newblock A new look at the classical sequence spaces by using multiplicative calculus.
\newblock {\em Politehn. Univ. Buchar. Sci. Bull. Ser. A Appl. Math. Phys.}
{\bf 2016}, {\em 78},~9--20.

\bibitem[Leitmann(1981)]{MR641031}
Leitmann, G.
\newblock {\em The Calculus of Variations and Optimal Control; 
Mathematical Concepts and Methods in Science and Engineering}; 
Plenum Press: New York, NY, USA; London, UK, 1981; p. xvi+311.

\bibitem[van Brunt(2004)]{MR2004181}
van Brunt, B.
\newblock {\em The Calculus of Variations}; Universitext, Springer: 
New York, NY, USA, 2004; p. xiv+290, doi:10.1007/b97436.

\bibitem[Sidi~Ammi and Torres(2008)]{MR2437854}
Sidi~Ammi, M.R.; Torres, D.F.M.
\newblock Regularity of solutions to higher-order 
integrals of the calculus of variations.
\newblock {\em Int. J. Syst. Sci.} {\bf 2008}, {\em 39},~889--895, 
doi:10.1080/00207720802184733.
\newblock {\tt arXiv:0707.2816}

\bibitem[Almeida \em{et~al.}(2019)Almeida, Tavares, and Torres]{MR3822307}
Almeida, R.; Tavares, D.; Torres, D.F.M.
\newblock {\em The Variable-Order Fractional Calculus of Variations}; 
Springer: Cham, Switzerland,  2019; pp.~xiv+124, 
doi:10.1007/978-3-319-94006-9.
\newblock {\tt arXiv:1805.00720}

\bibitem[Almeida and Martins(2020)]{MR4159537}
Almeida, R.; Martins, N.
\newblock Fractional variational principle of {H}erglotz for a new class of
problems with dependence on the boundaries and a real parameter.
\newblock {\em J. Math. Phys.} {\bf 2020}, {\em 61},~102701, 
doi:10.1063/5.0021373.

\bibitem[Malinowska and Torres(2014)]{MR3184533}
Malinowska, A.B.; Torres, D.F.M.
\newblock {\em Quantum Variational Calculus}; 
Springer: Cham, Switzerland, 2014, doi:10.1007/978-3-319-02747-0.

\bibitem[Brito~da Cruz and Martins(2018)]{MR3771533}
Brito~da Cruz, A.M.C.; Martins, N.
\newblock General quantum variational calculus.
\newblock {\em Stat. Optim. Inf. Comput.} {\bf 2018}, {\em 6},~22--41, 
doi:10.19139/soic.v6i1.467.

\bibitem[Martins and Torres(2009)]{MR2671876}
Martins, N.; Torres, D.F.M.
\newblock Calculus of variations on time scales with nabla derivatives.
\newblock {\em Nonlinear Anal.} {\bf 2009}, {\em 71},~e763--e773, 
doi:10.1016/j.na.2008.11.035.
\newblock {\tt arXiv:0807.2596}

\bibitem[Dryl and Torres(2017)]{MR3718404}
Dryl, M.; Torres, D.F.M.
\newblock Direct and inverse variational problems on time scales: A survey. 
In {\em Modeling, Dynamics, Optimization and Bioeconomics, {II}, 
Proceedings of the International Conference on Dynamics, Games 
and Science, Porto, Portugal, 17--21 February 2014}; 
Springer Proceedings in Mathematics \& Statistics; 
Springer: Cham, Switzerland, 2017; Volume 195, pp. 223--265, 
doi:10.1007/978-3-319-55236-1\_12.
\newblock {\tt arXiv:1601.05111}

\bibitem[Silva and Torres(2006)]{MR2323273}
Silva, C.J.; Torres, D.F.M.
\newblock Two-dimensional {N}ewton's problem of minimal resistance.
\newblock {\em Control Cybernet.} {\bf 2006}, {\em 35},~965--975.
\newblock {\tt arXiv:math/0607197}

\bibitem[Frederico \em{et~al.}(2014)Frederico, Odzijewicz, and Torres]{MR3169197}
Frederico, G.S.F.; Odzijewicz, T.; Torres, D.F.M.
\newblock Noether's theorem for non-smooth extremals 
of variational problems with time delay.
\newblock {\em Appl. Anal.} {\bf 2014}, {\em 93},~153--170, 
doi:10.1080/00036811.2012.762090.
\newblock {\tt arXiv:1212.4932}

\bibitem[Ferreira and Torres(2007)]{MR2405376}
Ferreira, R.A.C.; Torres, D.F.M.
\newblock Remarks on the calculus of variations on time scales.
\newblock {\em Int. J. Ecol. Econ. Stat.} {\bf 2007}, {\em 9},~65--73.
\newblock {\tt arXiv:0706.3152}

\bibitem[Dryl and Torres(2014)]{MR3313737}
Dryl, M.; Torres, D.F.M.
\newblock A general delta-nabla calculus of variations on time scales 
with application to economics.
\newblock {\em Int. J. Dyn. Syst. Differ. Equ.} {\bf 2014}, {\em 5},~42--71, 
doi:10.1504/IJDSDE.2014.067108.
\newblock {\tt arXiv:1410.1190}

\bibitem[Jost(2014)]{MR3157168}
Jost, J.
\newblock {\em Mathematical Methods in Biology and Neurobiology}; 
Universitext, Springer: London, UK, 2014; 
p. x+226, doi:10.1007/978-1-4471-6353-4.

\bibitem[Lemos-Pai\~{a}o \em{et~al.}(2020)Lemos-Pai\~{a}o, Silva, Torres, and
Venturino]{MR4110649}
Lemos-Pai\~{a}o, A.P.; Silva, C.J.; Torres, D.F.M.; Venturino, E.
\newblock Optimal control of aquatic diseases: 
A case study of {Y}emen's cholera outbreak.
\newblock {\em J. Optim. Theory Appl.} {\bf 2020}, {\em 185},~1008--1030, 
doi:10.1007/s10957-020-01668-z.
\newblock {\tt arXiv:2004.07402}

\bibitem[Grbi\'{c} \em{et~al.}(2017)Grbi\'{c}, Medi\'{c}, Perovi\'{c},
Mihailovi\'{c}, Novkovi\'{c}, and Durakovi\'{c}]{MR3609374}
Grbi\'{c}, T.; Medi\'{c}, S.; Perovi\'{c}, A.; Mihailovi\'{c}, B.;
Novkovi\'{c}, N.; Durakovi\'{c}, N.
\newblock A premium principle based on the {$g$}-integral.
\newblock {\em Stoch. Anal. Appl.} {\bf 2017}, {\em 35},~465--477, 
doi:10.1080/07362994.2016.1267574.

\bibitem[\"{U}nl\"{u}yol \em{et~al.}(2017)\"{U}nl\"{u}yol, Sala\c{s}, 
and \.{I}\c{s}can]{MR3750265}
\"{U}nl\"{u}yol, E.; Sala\c{s}, S.; \.{I}\c{s}can, I.
\newblock A new view of some operators and their properties 
in terms of the non-{N}ewtonian calculus.
\newblock {\em Topol. Algebra Appl.} {\bf 2017}, {\em 5},~49--54, 
doi:10.1515/taa-2017-0008.

\bibitem[Tekin and Ba\c{s}ar(2013)]{MR3073475}
Tekin, S.; Ba\c{s}ar, F.
\newblock Certain sequence spaces over the non-{N}ewtonian complex field.
\newblock {\em Abstr. Appl. Anal.} {\bf 2013}, \emph{2013}, 
739319, doi:10.1155/2013/739319.

\bibitem[Grossman(1979)]{MR557734}
Grossman, M.
\newblock {\em The First Nonlinear System of Differential 
and Integral Calculus}; MATHCO: Rockport, MA, USA,  1979; p. xi+85.

\bibitem[Grossman(1983)]{MR695495}
Grossman, M.
\newblock {\em Bigeometric Calculus}; Archimedes Foundation: 
Rockport, MA, USA, 1983; p. vii+100.

\bibitem[Riza and Akt\"{o}re(2015)]{MR3389882}
Riza, M.; Akt\"{o}re, H.
\newblock The {R}unge-{K}utta method in geometric multiplicative calculus.
\newblock {\em LMS J. Comput. Math.} {\bf 2015}, {\em 18},~539--554, 
doi:10.1112/S1461157015000145.

\bibitem[Boruah and Hazarika(2021)]{MR4222466}
Boruah, K.; Hazarika, B.
\newblock Some basic properties of bigeometric calculus 
and its applications in numerical analysis.
\newblock {\em Afr. Mat.} {\bf 2021}, {\em 32},~211--227, 
doi:10.1007/s13370-020-00821-1.

\bibitem[Bashirov and Norozpour(2017)]{MR3661747}
Bashirov, A.E.; Norozpour, S.
\newblock On complex multiplicative integration.
\newblock {\em TWMS J. Appl. Eng. Math.} {\bf 2017}, {\em 7},~82--93.

\bibitem[Waseem \em{et~al.}(2018)Waseem, Aslam~Noor, 
Ahmed~Shah, and Inayat~Noor]{MR3794497}
Waseem, M.; Aslam~Noor, M.; Ahmed~Shah, F.; Inayat~Noor, K.
\newblock An efficient technique to solve nonlinear equations using
multiplicative calculus.
\newblock {\em Turk. J. Math.} {\bf 2018}, {\em 42},~679--691, 
doi:10.3906/mat-1611-95.

\bibitem[Ali \em{et~al.}(2019)Ali, Abbas, and Zafar]{MR4016182}
Ali, M.A.; Abbas, M.; Zafar, A.A.
\newblock On some {H}ermite-{H}adamard integral inequalities in multiplicative calculus.
\newblock {\em J. Inequal. Spec. Funct.} {\bf 2019}, {\em 10},~111--122.

\bibitem[C\'{o}rdova-Lepe(2009)]{MR2724186}
C\'{o}rdova-Lepe, F.
\newblock From quotient operation toward a proportional calculus.
\newblock {\em Int. J. Math. Game Theory Algebra} {\bf 2009}, 
{\em 18},~527--536.

\bibitem[Yaz\i~c\i and Selvitopi(2017)]{MR3725737}
Yaz\i c\i, M.; Selvitopi, H.
\newblock Numerical methods for the multiplicative partial differential equations.
\newblock {\em Open Math.} {\bf 2017}, {\em 15},~1344--1350, 
doi:10.1515/math-2017-0113.

\bibitem[Boruah and Hazarika(2018)]{MR3813861}
Boruah, K.; Hazarika, B.
\newblock {$G$}-calculus.
\newblock {\em TWMS J. Appl. Eng. Math.} {\bf 2018}, {\em 8},~94--105.

\bibitem[Pap(2008)]{PAP2008368}
Pap, E.
\newblock Generalized real analysis and its applications.
\newblock {\em Int. J. Approx. Reason.} {\bf 2008}, 
{\em 47},~368--386, doi:10.1016/j.ijar.2007.05.015.

\bibitem[Duyar \em{et~al.}(2015)Duyar, Sa{\u{g}}Ir, 
and O{\u{g}}ur]{duyar2015some}
Duyar, C.; Sa{\u{g}}Ir, B.; O{\u{g}}ur, O.
\newblock Some basic topological properties on {N}on-{N}ewtonian real line.
\newblock {\em Br. J. Math. Comput. Sci.} {\bf 2015}, {\em 9},~300--307.

\bibitem[Binba\c{s}\i~o\v{g}lu \em{et~al.}(2016)Binba\c{s}\i~o\v{g}lu, 
Demiriz, and T\"{u}rko\v{g}lu]{MR3463533}
Binba\c{s}\i o\v{g}lu, D.; Demiriz, S.; T\"{u}rko\v{g}lu, D.
\newblock Fixed points of non-{N}ewtonian contraction mappings on
non-{N}ewtonian metric spaces.
\newblock {\em J. Fixed Point Theory Appl.} {\bf 2016}, 
{\em 18},~213--224, doi:10.1007/s11784-015-0271-y.

\bibitem[G\"{u}ng\"{o}r(2020)]{MR4104356}
G\"{u}ng\"{o}r, N.
\newblock Some geometric properties of the non-{N}ewtonian sequence spaces {$l_p(N)$}.
\newblock {\em Math. Slovaca} {\bf 2020}, {\em 70},~689--696, doi:10.1515/ms-2017-0382.

\bibitem[Sa\u{g}\i~r and Erdo\u{g}an(2020)]{MR4170025}
Sa\u{g}\i r, B.; Erdo\u{g}an, F.
\newblock On non-{N}ewtonian power series and its applications.
\newblock {\em Konuralp J. Math.} {\bf 2020}, {\em 8},~294--303.

\bibitem[Burgin and Czachor(2021)]{book2021:Burgin:Czachor}
Burgin, M.; Czachor, M.
\newblock {\em Non-diophantine Arithmetics in Mathematics, 
Physics and Psychology}; World Scientific: Singapore, 2021.

\bibitem[Plakhov and Torres(2005)]{MR2164553}
Plakhov, A.Y.; Torres, D.F.M.
\newblock Newton's aerodynamic problem in media of chaotically moving particles.
\newblock {\em Mat. Sb.} {\bf 2005}, {\em 196},~111--160, 
doi:10.1070/SM2005v196n06ABEH000904.
\newblock {\tt arXiv:math/0407406}

\bibitem[Torres and Plakhov(2006)]{MR2201075}
Torres, D.F.M.; Plakhov, A.Y.
\newblock Optimal control of {N}ewton-type problems of minimal resistance.
\newblock {\em Rend. Semin. Mat. Univ. Politec. Torino} {\bf 2006}, 
{\em 64},~79--95.
\newblock {\tt arXiv:math/0404237}

\bibitem[Castillo \em{et~al.}(2008)Castillo, Luce\~{n}o, and Pedregal]{MR2378083}
Castillo, E.; Luce\~{n}o, A.; Pedregal, P.
\newblock Composition functionals in calculus of variations. 
{A}pplication to products and quotients.
\newblock {\em Math. Model. Methods Appl. Sci.} {\bf 2008}, 
{\em 18},~47--75, doi:10.1142/S0218202508002607.

\bibitem[Riza \em{et~al.}(2009)Riza, \"{O}zyapici, and Misirli]{MR2588234}
Riza, M.; \"{O}zyapici, A.; Misirli, E.
\newblock Multiplicative finite difference methods.
\newblock {\em Quart. Appl. Math.} {\bf 2009}, {\em 67},~745--754, 
doi:10.1090/S0033-569X-09-01158-2.

\bibitem[Yener and Emiroglu(2015)]{MR3461691}
Yener, G.; Emiroglu, I.
\newblock A {$q$}-analogue of the multiplicative calculus: 
{$q$}-multiplicative calculus.
\newblock {\em Discret. Contin. Dyn. Syst. Ser. S} {\bf 2015}, 
{\em 8},~1435--1450, doi:10.3934/dcdss.2015.8.1435.

\bibitem[Stanley(1999)]{stanley1999multiplicative}
Stanley, D.
\newblock A multiplicative calculus.
\newblock {\em Primus} {\bf 1999}, {\em 9},~310--326.

\bibitem[Filip and Piatecki(2014)]{MR3176131}
Filip, D.A.; Piatecki, C.
\newblock A non-newtonian examination 
of the theory of exogenous economic growth.
\newblock {\em Math. Aeterna} {\bf 2014}, {\em 4},~101--117.

\bibitem[Uzer(2010)]{MR2734313}
Uzer, A.
\newblock Multiplicative type complex calculus as an alternative to the classical calculus.
\newblock {\em Comput. Math. Appl.} {\bf 2010}, {\em 60},~2725--2737, 
doi:10.1016/j.camwa.2010.08.089.

\bibitem[\c{C}akmak and Ba\c{s}ar(2014)]{MR3198166}
\c{C}akmak, A.F.; Ba\c{s}ar, F.
\newblock Certain spaces of functions over the field of non-{N}ewtonian complex numbers.
\newblock {\em Abstr. Appl. Anal.} {\bf 2014}, \emph{2014},  236124, doi:10.1155/2014/236124.

\bibitem[Bashirov and Norozpour(2018)]{MR3891225}
Bashirov, A.E.; Norozpour, S.
\newblock On an alternative view to complex calculus.
\newblock {\em Math. Methods Appl. Sci.} {\bf 2018}, 
{\em 41},~7313--7324, doi:10.1002/mma.4827.

\bibitem[\c{C}akmak and Ba\c{s}ar(2012)]{MR3016664}
\c{C}akmak, A.F.; Ba\c{s}ar, F.
\newblock Some new results on sequence spaces with respect to non-{N}ewtonian calculus.
\newblock {\em J. Inequal. Appl.} {\bf 2012}, \emph{2012}, 228, 
doi:10.1186/1029-242X-2012-228.

\bibitem[Florack and van Assen(2012)]{MR2878647}
Florack, L.; van Assen, H.
\newblock Multiplicative calculus in biomedical image analysis.
\newblock {\em J. Math. Imaging Vis.} {\bf 2012}, {\em 42},~64--75, 
doi:10.1007/s10851-011-0275-1.

\end{thebibliography}
\end{document}